\documentclass{amsart}

\usepackage{times, epsfig, amsmath, amsthm, amssymb, mathrsfs}

\begin{document}

\newtheorem{thm}{Theorem}[subsection]
\newtheorem{lem}[thm]{Lemma}
\newtheorem{cor}[thm]{Corollary}
\newtheorem{conj}[thm]{Conjecture}
\newtheorem{qn}[thm]{Question}
\newtheorem*{no_smooth_action_thm}{Theorem A}
\newtheorem*{representation_topology_cor}{Corollary D}
\newtheorem*{transcendental_module_cor}{Corollary B}
\newtheorem*{outer_approximation_thm}{Theorem C}

\theoremstyle{definition}
\newtheorem{defn}[thm]{Definition}
\newtheorem{note}[thm]{Notation}

\theoremstyle{remark}
\newtheorem{rmk}[thm]{Remark}
\newtheorem{exa}[thm]{Example}

\def\R{\mathbb R}
\def\C{\mathbb C}
\def\Z{\mathbb Z}
\def\N{\mathbb N}
\def\Q{\mathbb Q}
\def\H{\mathbb H}
\def\CP{\mathbb{CP}}
\def\RP{\mathbb{RP}}
\def\O{\mathscr{O}}
\def\id{{\text{Id}}}
\def\tr{{\text{tr}}}
\def\hom{{\text{Hom}}}
\def\Gal{{\text{Gal}}}
\def\homeo{{\text{Homeo}}}
\def\rot{{\text{rot}}}
\def\SL{{\text{SL}}}
\def\SU{{\text{SU}}}
\def\PSL{{\text{PSL}}}
\def\GL{{\text{GL}}}
\def\til{\widetilde}

\title{Dynamical forcing of circular groups}

\author{Danny Calegari}

\address{Department of Mathematics \\ California Institute of Technology
\\ Pasadena CA 91125}
\email{dannyc@its.caltech.edu}
\date{8/31/2004; version 0.09} 

\begin{abstract}
In this paper we introduce and study the notion of {\em dynamical forcing}.
Basically, we develop a toolkit of techniques to produce finitely
presented groups which can only act on the circle with
certain prescribed dynamical properties.

As an application, we show that the set $X \subset \R/\Z$ consisting
of rotation numbers $\theta$ which can be {\em forced} by
finitely presented groups is an infinitely generated $\Q$--module,
containing countably infinitely
many algebraically independent transcendental numbers.
Here a rotation number $\theta$ is {\em forced} by a pair
$(G_\theta,\alpha)$ where $G_\theta$ is a finitely
presented group $G_\theta$ and $\alpha \in G_\theta$ is some element,
if the set of rotation numbers of $\rho(\alpha)$ as $\rho$
varies over $\rho \in \hom(G_\theta,\homeo^+(S^1))$ is
precisely the set $\lbrace 0, \pm \theta \rbrace$.

We show that the set of subsets of $\R/\Z$ which are of the form
$$\rot(X(G,\alpha)) = \lbrace r \in \R/\Z \; | \; r = \rot(\rho(\alpha)),
\rho \in \hom(G,\homeo^+(S^1)) \rbrace$$
where $G$ varies over countable groups, are exactly the set of
closed subsets which contain $0$ and are invariant under $x \to -x$.
Moreover, we show that every such subset can be approximated from
above by $\rot(X(G_i,\alpha_i))$ for {\em finitely presented} $G_i$.

As another application, we construct a finitely generated group $\Gamma$ which
acts faithfully on the circle, but which does not admit any faithful $C^1$ action,
thus answering in the negative a question of John Franks.
\end{abstract}

\maketitle

\section{Introduction}

It is a basic problem, given a group $G$, to understand and classify actions
of $G$ on a topological space $X$. For good reasons, one of the
most widely studied cases is when $G$ is an explicit finitely presented
group or a representative of a class of groups, and $X$ is the circle $S^1$.
One then tries to understand what possible dynamics $G$, or some element
of $G$, can have under representations $G \to \homeo^+(S^1)$.

In this paper, we treat a complementary problem --- given a collection of dynamical
constraints, we study when we can produce a finitely presented group $G$
whose representations to $\homeo^+(S^1)$ are subject to precisely this
collection of constraints. Explicitly, we are interested in characterizing which
sets of rotation numbers can be {\em forced} by finitely presented groups.
That is, given a group $G$ and an element $\alpha \in G$, what possible subsets
of the circle can arise as the set of rotation numbers of $\alpha$ as we vary
over all representations of $G$ into $\homeo^+(S^1)$.

This is part of a broader program to develop a {\em toolkit} for producing
finitely generated groups with prescribed algebraic and geometric properties,
and for embedding given groups in larger groups in which desirable properties
of the smaller group persist, but undesirable ones do not.

\subsection{Statement of results}

In section \S 2 and \S 3, we warm up by describing a finitely generated group
which acts faithfully on $S^1$, but admits no faithful $C^1$ action. That is,

\begin{no_smooth_action_thm}
There exists a finitely generated group which acts faithfully on $S^1$
by orientation preserving homeomorphisms, but for which no faithful action is
conjugate to a $C^1$ action.
\end{no_smooth_action_thm}

The group $Q$ in question contains as a subgroup $\Gamma$, the
fundamental group of the unit tangent bundle of the hyperbolic $(2,3,7)$
orbifold, a well--known perfect $3$--manifold group. $Q$ is a quotient
of $\widehat{\Gamma}$, which is obtained from $\Gamma$ by adding $3$ generators
and $3$ relations. The existence of this example answers in the negative 
a question posed by John Franks.

In section \S 4 we initiate a systematic study of the possible sets of rotation
numbers which can be {\em forced} by algebraic properties of a group.

Explicitly, let $G$ be a countable group, and $\alpha$ an element of $G$.
We study the set of values 
$$\rot(X(G,\alpha)) = \bigcup_\rho \rot(\rho(\alpha)) \subset \R/\Z$$
where $\rho$ varies over all representations $\rho:G \to \homeo^+(S^1)$,
and $\rot$ denotes the {\em rotation number} of an element of $\homeo^+(S^1)$,
as defined by Poincar\'e. We determine exactly which subsets
$\rot(X(G,\alpha))$ can occur, and show that they define a topology on
$\R/\Z$ called the {\em representation topology}.

We prove the following result:

\begin{representation_topology_cor}
The set of subsets of $S^1$ of the form $\rot(X(G,\alpha))$ 
where $G$ varies over all countable groups, and $\alpha \in G$ is
arbitrary, are precisely the nonempty closed subsets of a topology, called the
{\em representation topology}.

The nonempty closed subsets in the representation topology on $S^1$ are
exactly unions $\lbrace 0 \rbrace \cup K$ where $K$ is closed (in the
usual sense) and invariant under $x \to -x$.
\end{representation_topology_cor}

More interesting (and much more difficult) is to characterize which
subsets of $\R/\Z$ can arise as rotation numbers of some element $\alpha$ of
{\em finitely presented} groups $G$. Say that a number $\theta \in \R/\Z$ is
{\em forceable} if $\rot(X(G,\alpha)) = \lbrace 0 \cup \pm \theta \rbrace$
for some finitely presented group $G$. Notice that since there are only
countably many finitely presented groups, there are only countably many
forceable numbers.

To describe our results, we must first introduce the
operation $+_l$ for $l \in \R$. This is a (partially defined)
symmetric binary operation on $\R/\Z$. A quick way to define it is as follows:
if $\alpha_1,\alpha_2$ are rotations through (signed)
angles $2\pi\theta_1,2\pi\theta_2$
of the hyperbolic plane $\H^2$ with
centers at points $p_1,p_2$ which are distance $l$ apart, then
when the composition $\alpha_1 \circ \alpha_2$ is a rotation, it is a rotation
through angle $2\pi (\theta_1 +_l \theta_2)$. 

If $Y \subset S^1$ is a set of rotation numbers and $L \subset \R$, 
the {\em algebraic closure of $Y$} with respect to the operations $+_l$
with $l \in L$, is defined inductively as
the smallest set which includes $Y$ itself, and also includes every
solution to every {\em almost determined} finite system of equations
with operations $+_l$, and coefficients in the algebraic closure of
$Y$. Here a system of equations is {\em almost determined} if it has only
finitely many solutions.

With this definition,
our main results are the following:

\begin{transcendental_module_cor}
Let $X$ be the set of forceable rotation numbers. Then $X$ contains
countably infinitely many algebraically independent transcendental
numbers, as well as all the rational numbers. Moreover,
there is a dense set $L$ of real numbers $l \in \R$ containing $0$,
which are of the form
$\log(r)$ for $r$ algebraic, so that $X$ is algebraically closed
with respect to the operations $+_l$.
\end{transcendental_module_cor}

In particular, the set of forceable rotation numbers is the
reduction mod $\Z$ of an infinite dimensional $\Q$--vector subspace
of $\R$, which generates a field of infinite transcendence degree over $\Q$.

Finally, we show that every closed subset in the representation topology
can be approximated (in a constructive way) by closed subsets forced by
finitely presented groups.

\begin{outer_approximation_thm}
Let  $K$ be any closed subset of $\R/\Z$. Then there are
a sequence of pairs $G_{K_i},\alpha_i$ where $G_{K_i}$ is
a finitely presented group and $\alpha_i \in G_{K_i}$, and
closed subsets $K_i$ of $\R/\Z$ such that
\begin{enumerate}
\item{Each $K_{i+1} \subset K_i$.}
\item{The intersection $\cap_i K_i = K$.}
\item{The set of rotation numbers $\rho(\alpha_i)$ as $\rho$ ranges over
$\hom(G_{K_i},\homeo^+(S^1))$ is exactly equal to 
$\lbrace 0 \rbrace \cup K_i \cup - K_i$.}
\item{There is a canonical element
$\nu_i \in G_{K_i}$ of order $3$, so that if $\rot(\rho(\nu_i))=0$, $\rho(\alpha_i)$
is trivial, and if $\rot(\rho(\nu_i))=1/3$, the set
of compatible $\rot(\rho(\alpha_i))$ is exactly equal to $K_i$.}
\end{enumerate}
\end{outer_approximation_thm}

Note that the canonical element $\nu_i$ can be thought of as a ``dial'' which
can be set to three possible values, which resolve the necessary ambiguity
in the set $\rot(X(G_{K_i},\alpha_i))$, namely that it be invariant under
$x \to -x$, and that it contain $0$.

\subsection{Acknowledgements}

The subject matter in this paper was partly inspired by a discussion with John
Franks and Amie Wilkinson, and a comment in an email from \'Etienne
Ghys. I would like to thank Nathan Dunfield for some excellent comments
on an earlier version of this paper, and Hee Oh for some useful info about arithmetic
lattices. Thanks as well to the referee, for catching
a number of errors, especially in some of the formulae.

\section{A non--smoothable group action on $S^1$}

Let $\Gamma$ be a finitely generated group of homeomorphisms of the circle.
Then $\Gamma$ quasi--preserves a harmonic measure with support contained in
any invariant set. If the action of $\Gamma$ on $S^1$ is minimal, by integrating
this measure, one sees that $\Gamma$ is conjugate to a Lipschitz action.

Not every action is conjugate to a $C^1$ action; however, John Franks asked
whether every abstract group $\Gamma$ of homeomorphisms of $S^1$ was
abstractly isomorphic to a group of $C^1$ diffeomorphisms. We will answer this
question in the negative.

The first step is to find a particular group $\Gamma$ and a faithful
{\em action} of this group on $S^1$ which is not conjugate to a $C^1$ action.
Then we embed the group $\Gamma$ in a larger group $\widehat{\Gamma}$,
adding finitely many
generators and relations, which force every action of the group $\widehat{\Gamma}$ 
on $S^1$ to restrict to an action by $\Gamma$ of this nonsmoothable kind.

Let $\Gamma$ be the fundamental group of the unit tangent bundle of
the hyperbolic $(2,3,7)$--orbifold. A presentation for $\Gamma$ is given by:

$$\Gamma = \langle A,B,C,T \; | \; A^2 = T, B^3 = T, C^7 = T, ABC = T \rangle$$

Then $T$ generates a free cyclic central subgroup. Quotienting by $\langle T \rangle$
gives a homomorphism from $\Gamma$ to the fundamental group of the $(2,3,7)$ triangle
orbifold $O$, which we denote by $\Delta$. A presentation for $\Delta$ is

$$\Delta = \langle A,B,C \; | \; A^2 = B^3 = C^7 = ABC = \id \rangle$$

Obviously, these presentations are not minimal.
Note that $\Gamma$ is a {\em perfect group}; that is, $H_1(\Gamma;\Z) = 0$.

So we have a short exact sequence
$$ 0 \to \Z \to \Gamma \to \Delta \to 0$$

The group $\Delta$ admits a faithful representation in $\PSL(2,\R)$, coming
from its realization as the fundamental group of a hyperbolic $2$--orbifold.
This real analytic action can be Denjoyed, by blowing up an orbit with
trivial stabilizer, and inserting an action of $\Z$, in a $C^1$ manner, 
to give a faithful $C^1$ action of $\Gamma$ on $S^1$. Note that since 
every nontrivial element of $\PSL(2,\R)$ has at most two fixed points, all
but countably many points in $S^1$ have trivial stabilizer.
See \cite{Den} for Denjoy's construction.

In another way, we can think of $\Gamma$ as the universal central extension of
$\Delta$, and identify it with the preimage of $\Delta$ in the universal
covering group $\til{\PSL(2,\R)}$ of $\PSL(2,\R)$. In this way, we get
a faithful real analytic action of $\Gamma$ on $\R$. By inserting $\R$ into
$S^1$, we can produce many faithful homomorphisms from $\Gamma$ to $\homeo^+(S^1)$.

The following lemma gives a simple condition under which such a homomorphism
from $\Gamma$ to $\homeo^+(S^1)$ must have a global fixed point.

\begin{lem}\label{common_fixed_point}
Let $h:\Gamma \to \homeo^+(S^1)$ be a faithful action
such that $A,B,C$ all have fixed points.
Then $\Gamma$ has a global fixed point.
\end{lem}
\begin{proof}
Let $\Lambda_A,\Lambda_B,\Lambda_C$ be the fixed point sets of $A,B,C$
respectively. Then each of $A,B,C$ is conjugate to a translation
on the complementary intervals of their fixed point sets, and the same
is true of their positive powers. In particular, since $A,B,C$ have
common positive powers, the fixed point sets of $A,B,C$ are equal to the
fixed point set of $T$ and to each other. It follows that
$$\Lambda_A = \Lambda_B = \Lambda_C = \text{fix}(\Gamma)$$
is nonempty.
\end{proof}

A group is {\em locally indicable} if every nontrivial finitely generated subgroup
admits a surjective homomorphism to $\Z$.

The following theorem is known as the {\em Thurston stability theorem} \cite{Thur}:

\begin{thm}[Thurston]
Let $G$ be a group of germs at $0$ of orientation preserving
$C^1$ diffeomorphisms of $\R$. Then $G$ is locally indicable.
\end{thm}

It follows that an action of $\Gamma$ on $S^1$ for which all $A,B,C$
have fixed points is not conjugate to a $C^1$ action.

\section{Forcing dynamics of $\Gamma$}

We will show how to embed $\Gamma$ in a larger finitely generated group
which still acts faithfully on $S^1$, such that every action of
the larger group restricts to an action of $\Gamma$ with a global fixed point.

Define the following group

$$\widehat{\Gamma} = \langle \Gamma, X,Y,Z \; | 
\; XAX^{-1} = A^2, YBY^{-1} = B^2, ZCZ^{-1} = C^2 \rangle$$

Notice that $\widehat{\Gamma}$ contains a number of copies of
the Baumslag--Solitar group $BS_2$.

In the following lemma, we show that every orientation--preserving
action of $\widehat{\Gamma}$ on $S^1$ restricts on $\Gamma$ to an action
satisfying the properties of lemma~\ref{common_fixed_point}. But first we
recall Poincar\'e's definition of a {\em rotation number}.

\begin{defn}[Poincar\'e]
Let $\alpha \in \homeo^+(S^1)$, and let $\til{\alpha} \in \homeo^+(\R)$ be
a lift of $\alpha$. Define the rotation number $\rot(\alpha)$ by
$$\rot(\alpha) = \lim_{n \to \infty} \frac {\til{\alpha}^n(0)} {n} \mod \Z$$
\end{defn}

It is easy to deduce directly from the definition that
the rotation number is a continuous class function from $\homeo^+(S^1) \to \R/\Z$.
Moreover the formula $\rot(\alpha^n) = n \, \rot(\alpha)$ holds
for any integer $n$, and $\rot(\alpha) = 0$ iff $\alpha$ has a fixed point.
A basic reference is \cite{Poin}.

\begin{lem}\label{rotation}
For every homomorphism of $\widehat{\Gamma}$ to $\homeo^+(S^1)$, the
elements $A,B,C$ have a fixed point.
\end{lem}
\begin{proof}
By the properties of Poincar\'e's rotation number, we calculate
$$\rot(A) = \rot(XAX^{-1}) = \rot(A^2) = 2\rot(A)$$
and similarly for $B$ and $C$. In particular, the rotation numbers of $A,B,C$
are zero, and therefore they all have fixed points.
\end{proof}

\begin{no_smooth_action_thm}
There exists a finitely generated group which acts faithfully on $S^1$
by orientation preserving homeomorphisms, but for which no faithful action is
conjugate to a $C^1$ action.
\end{no_smooth_action_thm}
\begin{proof}
Let $h: \Gamma \to \homeo^+(S^1)$ be a faithful action obtained by compactifying
an action on $\R$ arising from a faithful representation into $\til{\PSL(2,\R)}$.
Since each $A,B,C$ acts on $S^1$ with a single fixed point, each is conjugate
to any positive power of itself. It follows that $h$ extends to a homomorphism
$h:\widehat{\Gamma} \to \homeo^+(S^1)$. Let $Q$ be the image of $\widehat{\Gamma}$
under this homomorphism.

By lemma~\ref{rotation}, for any faithful representation of $Q$ in $\homeo^+(S^1)$,
the elements $A,B,C$ must have fixed points. It follows by
lemma~\ref{common_fixed_point} that $\Gamma$ has a global common fixed point.
By the Thurston stability theorem, the action of $\Gamma$ is not conjugate to a
$C^1$ action near a common fixed point in the frontier of $\text{fix}(\Gamma)$.
\end{proof}

\section{Forcing rotation numbers}

\subsection{Galois groups and the Milnor--Wood inequality}

We think of the group $\PSL(2,\R)$ as the group of real
projective automorphisms of the
circle $\RP^1$ and $\PSL(2,\C)$ as the group of complex projective
automorphisms of the (Riemann) sphere $\CP^1$.

The conjugacy class of a generic element $[\alpha] \in \PSL(2,\C)$ 
is determined by the square of the trace of a representative $\alpha$
unless this trace is $\pm 2$, in which case there are two possibilities.
Said another way, the map $$\tr^2:\PSL(2,\C) \to \C$$
is $1$--$1$ on conjugacy classes away from the point $4$.

For a finitely generated group $G$ and an element $\alpha \in G$, 
let $V(G) = \hom(G,\PSL(2,\C))$, and $V(\alpha) \subset \PSL(2,\C)$
be the union of the images of $\alpha$ under all the $\rho \in V(G)$.
Then define $X(\alpha) = \tr^2 V(\alpha)$. Here we emphasize
that we look at such sets for {\em all} pairs $(G,\alpha)$, and
not just different $\alpha$ contained in a fixed $G$.

Now, given any
$\rho:G \to \PSL(2,\C)$ such that $t = \tr^2(\rho(\alpha))$ and any
$\sigma \in \Gal(\C/\Q)$, we can look at $\rho^\sigma:G \to \PSL(2,\C)$,
and observe that 
$$\sigma(t) = \tr^2(\rho^\sigma(\alpha))$$
In particular, the subsets of $\C$ of the form
$X(\alpha)$ consist of all of $\C$, the
empty set, certain finite unions of Galois orbits of algebraic $k \in \C$, and
complements of certain finite unions of Galois orbits of algebraic $k \in \C$.
The reason complements of finite unions also arise is that although
the representation variety $V(\alpha)$ is an affine variety,
the morphism $\tr^2$ is not {\em proper},
and therefore the image is not always Zariski
closed, and may omit a (finite) Galois orbit in the closure.

Notice that by the Noetherian property of $\C[x]$ that the restriction
of the construction above to {\em finitely presented groups} determines the
same set of $X(\alpha)$.

In particular, the algebraic structure of $G$ constrains only
weakly the (topological) dynamics of a typical element $\alpha$ on
$\CP^1$.

From our point of view, the key difference between $\PSL(2,\R)$ and $\PSL(2,\C)$ is the
following inequality, known as the {\em Milnor--Wood inequality}:

\begin{thm}[Milnor--Wood inequality]\label{MilnorWood}
Let $E$ be a foliated circle bundle over a surface $\Sigma$ of
genus $g$. Then the Euler class of $E$ evaluated on the fundamental
class $[\Sigma]$ of $\Sigma$ satisfies
$$|e(E)([\Sigma])| \le \min (0, -\chi(\Sigma))$$ 
\end{thm}

By abuse of notation, we abbreviate the expression $e(E)([\Sigma])$ as
$e(E)$ and refer to it as the {\em Euler number} of the bundle $E$.

\subsection{Semi--conjugacy and Ghys' theorem}

\begin{defn}
A {\em monotone} map $S^1 \to S^1$ is a degree one map with connected
point inverses. Two group actions $\rho_i:G \to \homeo^+(S^1)$ for
$i=1,2$ are {\em semi--conjugate} if there is a third group action
$\rho:G \to \homeo^+(S^1)$ and monotone maps $\phi_i:S^1 \to S^1$
which intertwine the various group actions, so that $(\phi_i)_* \rho = \rho_i$
for $i=1,2$.
\end{defn}

It is not obvious from the definition that semi--conjugacy is an equivalence
relation. But it is implied by the existence of a {\em pushout diagram} for
every pair of monotone maps, which is natural, and therefore compatible with
$G$ actions (see \cite{CalDun} lemma 6.6 for a proof).

It is worth remarking that our notation is not standard: usually
one says $\rho_i,\rho_j$ are semi--conjugate if there is a monotone map
$\phi$ with $(\phi)_*\rho_i = \rho_j$. Semi--conjugacy as we have defined
it above is the equivalence relation generated by ordinary semi--conjugacy.
Some authors refer to semi--conjugacy as we have defined it as {\em monotone
equivalence}.

Note that the rotation number of an element of $G$ is invariant
under semi--conjugacy. In fact, it is not hard to show that rotation number is
a {\em complete} invariant of a homomorphism from $\Z$ to $\homeo^+(S^1)$ up
to semi--conjugacy.

Generalizing this, Ghys proved (\cite{Ghys},\cite{Matsumoto}) the following theorem:
\begin{thm}[Ghys]\label{semiconjugacy_invariant}
Let $G$ be a group. Then homomorphisms of $G$ to $\homeo^+(S^1)$ up to semi--conjugacy
are in bijective correspondence with elements of the second bounded
cohomology $[c] \in H^2_b(G;\Z)$ which admit representative cocycles $c$
taking values in $\lbrace 0,1 \rbrace$.
\end{thm}

Here the {\em bounded cohomology} of a group $G$ is the cohomology of the
complex of cochains $C_b^*$ which are bounded, thought of as functions on
the generators of $C_*$. In particular, the well--known fact that
$H^2_b(\Z;\Z) = S^1$ reflects the remark made above that rotation number determines
semiconjugacy on $\homeo^+(S^1)$.

For the group $\homeo^+(S^1)$ itself, it is known that
$$H^2_b(\homeo^+(S^1);\Z) = H^2(\homeo^+(S^1);\Z) = \Z$$
This follows from two facts: firstly, the cohomology of $\homeo^+(S^1)$
as a {\em discrete} group is equal to its cohomology as a {\em continuous}
group, by the fundamental theorem of Mather--Thurston (see \cite{Tsu} for
a nice exposition and references). Secondly, for any {\em uniformly
perfect} group $G$, $H^2_b \cong H^2$. Here, a group is said to be {\em uniformly
perfect} if there is some $n>0$ such that every element can be written as
a product of at most $n$ commutators. The group $\homeo^+(S^1)$ is easily
shown to be such a group.

The Milnor--Wood inequality implies that for a surface group
$\Gamma = \pi_1(\Sigma)$ where the genus of $\Sigma$ is at least $2$,
the maximal Euler class is equal to $-\chi(\Sigma)$. But the very
maximality of this Euler class forces the representative cocycle of an
action with maximal Euler class to be very constrained. From this, one
can easily deduce the following theorem of Matsumoto, proved in \cite{Matsumoto2}:

\begin{thm}[Matsumoto]\label{maximal_semiconjugate}
Let $\Sigma$ be an orientable surface of genus $\ge 2$, and let
$\Gamma = \pi_1(\Sigma)$. Let $\rho:\Gamma \to \homeo^+(S^1)$ be
an action of maximal Euler class. Then $\rho$ is semi--conjugate to a
Fuchsian action.
\end{thm}

We will make use of this theorem later.

For any group $G$, the trivial element of $H^2_b(G;\Z)$ corresponds to the
semi--conjugacy class of the trivial action of $G$ on $S^1$. 
Moreover, given any bounded cocycle $[c] \in H^2_b(G;\Z)$, 
the cocycle $[1-c]$ is homologous to $-[c]$. This reflects the fact that
any homomorphism $G \to \homeo^+(S^1)$ can be conjugated by an
{\em orientation--reversing} homeomorphism of $S^1$ to give another (typically
non--conjugate) homomorphism.

By analogy with \S 4.1 we define $X(S^1)$ to be the set of conjugacy classes of elements
of $\homeo^+(S^1)$. Then $\rot:X(S^1) \to S^1$ is well--defined.

\begin{defn}
For each pair $(G,\alpha)$ where
$G$ is a countable group, $\alpha \in G$ an element, let $X(G,\alpha)$ be the union of
the images in $X(S^1)$ of $\rho(\alpha)$ under all homomorphisms
$\rho:G \to \homeo^+(S^1)$. Then the {\em representation topology} on
$\R/\Z$ is the topology whose nonempty closed sets are {\em exactly} sets of the
form $\rot(X(G,\alpha))$ for $G,\alpha$ as above.
\end{defn}

It is not at all clear from the definition that this is really a topology ---
i.e. that this class of subsets of $S^1$ is closed under finite
union and arbitrary intersection. We will prove this in section \S 4.6.

Before we start, make the following definition:

\begin{defn}\label{dynamical_forcing_definition}
Given an abstract finitely presented group $\Gamma$ and a representation
$$\rho:\Gamma \to \homeo^+(S^1)$$
we say an embedding of $\Gamma$ in a finitely presented
group $G$ {\em forces $\rho$} if for every representation 
$$\sigma:G \to \homeo^+(S^1)$$
either $\sigma|_\Gamma$ is semi--conjugate to $\pm \rho$, or to $\id$. 
\end{defn}

We will be preoccupied with the case that $\Gamma = \Z$, generated by
a single element $\alpha$, and $\rho(\alpha)$ is a rotation through
angle $\theta$ for certain real numbers $\theta$, and we will say that
such a rotation number $\theta$ can be {\em forced}.

\subsection{Rational rotation numbers}

There are only countably many finitely presented groups, and therefore 
only countably many rotation numbers $\theta \in \R/\Z$ can
be forced. What are they? We enumerate several constructions which show that
the set of forceable rotation numbers has an algebraic structure which
is richer than one might initially expect.

At first glance, rational rotation numbers might seem easy to force.
In particular, for a rational number $r$ with reduced expression 
$r = p/q$, a rotation number with denominator $q$ can be forced by the relation
$$\alpha^q = \id$$
But of course, such an $\alpha$ could have rotation number $p/q$ for any $p$.
How then do we force a specific value of $p$? The cases $q=2,3,4$ are simple:
the possible choices of coprime $p$ are absorbed into the necessary ambiguity of
definition~\ref{dynamical_forcing_definition}. The next case we study is $q=7$,
skipping over $q=5$ for the moment. We want to find a group and an element which
can have rotation number $1/7$ for some action, but not $2/7$ or $3/7$.
 
Recall the group $\Delta$ from \S 2. $\Delta$ is the fundamental group of
the hyperbolic $(2,3,7)$--orbifold, and as such, it has an obvious
homomorphism to $\homeo^+(S^1)$ arising from its canonical representation in
$\PSL(2,\R) < \homeo^+(S^1)$. It will turn out that this representation is 
(up to semi--conjugacy and reflection) essentially the only nontrivial such
representation. But at this point we need only prove something weaker: that
the rotation number of the element of order $7$ must always be $1/7$
(up to sign) under any nontrivial homomorphism.

\begin{lem}\label{1/7}
If $\rho$ is any homomorphism from $\Delta$ to $\homeo^+(S^1)$, then either
$\rho$ is the trivial homomorphism, or else the rotation number of $\rho(C)$ is
equal to $1/7$.
\end{lem}
\begin{proof}
From the presentation of $\Delta$ from \S 2, we see that $\rho(A)$ must be
taken to $\pm \id$. If $\rho(A) = \id$, then $\rho(B) = \rho(C)$ because
$ABC^{-1}= \id$. Since the
orders of these elements is a factor of $3$ and $7$ respectively, they must
both be trivial. Similarly, if either of $\rho(C),\rho(B)$ is trivial, then
the whole group is. It follows that either $\rho$ is trivial, or, after
possibly conjugating by an orientation--reversing automorphism, we can assume the
rotation numbers of $\rho(A),\rho(B),\rho(C)$ are $1/2,1/3,p/7$
respectively for some $0 \le p \le 7$.

Now, the group $\Delta$ has an index 168 normal
subgroup $K$ which is the fundamental group of the Klein quartic $Q$, i.e.
a Riemann surface of genus $3$. This arises as a {\em congruence subgroup}
of $\Delta$; that is, there is a short exact sequence
$$0 \to K \to \Delta \to \PSL(2,\mathbb{F}_7) \to 0$$
where $\mathbb{F}_7$ denotes the field with $7$ elements.

The representation $\rho$ determines
a foliated bundle $E_O$ over the triangle orbifold $O$, which lifts to
a foliated bundle $E_Q$ over $Q$. We must be a bit careful about what we
mean here: a foliated $S^1$ bundle $E$ over a manifold $M$ is obtained from
a representation $\phi:\pi_1(M) \to \homeo^+(S^1)$ by taking a quotient
$$E = \til{M} \times S^1 / (m,\theta) \sim (\alpha(m),\phi(\alpha)(\theta))$$
for all $m \in \til{M},\theta \in S^1$ and $\alpha \in \pi_1(M)$.
The projection onto the first factor has fiber $S^1$ everywhere, and defines
a bundle structure. However, if $M$ is an {\em orbifold}, the action of the
orbifold fundamental group on the orbifold universal cover (if one exists) 
$\til{M}$ is not free, and $E$ is not quite a bundle in general. If $O$ is
a $2$--dimensional surface orbifold, then $E_O$ is actually a foliated
{\em Seifert fibered space} which fibers over the orbifold $O$. See
\cite{Montesinos} for more details. In any case, there is a well--defined
Euler number $|e(E_O)|$
for $E_O$ which is multiplicative under covers of the base.
The singular fibers --- corresponding to the orbifold points of $O$ --- contribute
fractional parts to the Euler number, equal to the rotation number of the
infinitesimal holonomy around an orbifold point. In particular,
$$|e(E_O)| = |n - 1/2 - 1/3 - p/7|$$
for some integer $n$. It follows that we can compute
$$|e(E_Q)| = 168|n - 1/2 - 1/3 - p/7|$$
By the Milnor--Wood inequality (theorem~\ref{MilnorWood}) $|e(E_Q)| \le 4$.
It follows that the only possibility for $p$ is $1$.
\end{proof}

As a consequence of this lemma, any rotation number $p/7$ can be forced by the group
$$\langle \Delta, \alpha \; | \; \alpha = C^p \rangle$$

Now we can use this example to force rotation numbers $p/q$ for any coprime $p,q$.

\begin{thm}\label{rational_rotation_numbers}
For any rational number $p/q$ with $p<q$ there is a finitely presented
group $G_{p/q}$ and an element $\alpha \in G_{p/q}$ such that
the set of rotation numbers of $\alpha$ under homomorphisms from $G_{p/q}$
to $\homeo^+(S^1)$ is exactly the set $\lbrace 0, p/q, -p/q \rbrace$.
\end{thm}
\begin{proof}
We show how to force $1/q$ for any $q>1$.
Define
$$S_q = \langle \mu, \nu, \alpha \; | \; [\mu,\nu] = \alpha, \alpha^q = \id \rangle$$
Then $S_q$ is the (orbifold) fundamental group of the genus $1$ orbifold $P$ with
one cone point of order $q$, and therefore 
contains a subgroup $S_q^*$ of index $2q$ which is isomorphic to
the fundamental group of $\Sigma_q$, the surface of genus $q$.
Then $\Sigma_q$ and the Klein quartic $Q$ have common covers which are
homeomorphic, and therefore have isomorphic fundamental groups $S_q^K, K^{S_q}$.
We can impose as generators the group $\Delta$ and an extra $\beta$, and relations
that $\beta K^{S_q} \beta^{-1} = S_q^K$.

Now, arguing as in the proof of lemma~\ref{1/7}, we know that
$$e(E_P) =\pm (n - p/q)$$ for some $p<q$, and therefore
$$e(E_{\Sigma_q}) = \pm 2q(n - p/q) \cong \pm -2p \mod 2q$$
On the other hand, the relation 
$\beta K^{S_q} \beta^{-1} = S_q^K$ fixes $|e(E_{\Sigma_q})|$,
setting it equal to exactly $2q-2$.
This implies $p=1$.

Note that this equality is realized for a representation of $S_q$ corresponding to
a hyperbolic orbifold structure on the torus with one cone point of order $q$;
this contains a surface group, whose action on $S^1$ is {\em topologically}
conjugate to the action of a suitable surface subgroup of $\Delta$. In particular,
the group we have constructed actually has an action on $S^1$ so that the
rotation number of $\alpha$ is $1/q$.

Of course once we can produce an element with rotation number exactly equal to
$0, 1/q, -1/q$ the $p$th power of such an element has rotation numbers as desired.
\end{proof}

Explicitly, a presentation for $G_{p/q}$ is

$$G_{p/q} = \langle \Delta, \mu,\nu, \beta, \gamma, \alpha \; | \;
[\mu,\nu] = \gamma, \beta K^{S_q} \beta^{-1} = S_q^K, \alpha = \gamma^p \rangle$$
which has $8$ generators, of which $2$ are clearly superfluous. 
Are these the ``simplest'' groups which force rational rotation numbers?

\begin{rmk}
In \cite{GhysSerg}, Ghys and Sergiescu show that Thompson's group admits a
$C^\infty$ action on $S^1$ which is
unique {\em up to $C^2$ conjugacy}. This group contains elements
of every possible rational rotation number. It seems likely that
this action is unique up to semi--conjugacy amongst all $C^0$ actions, but we
have not pursued this. 
\end{rmk}

\begin{rmk}
Notice by theorem~\ref{maximal_semiconjugate} that every action of
$\Delta$ is either trivial or semi--conjugate to a Fuchsian action, as remarked above.
\end{rmk}

\begin{rmk}
It is actually possible to force all rational rotation numbers with {\em torsion free}
finitely presented groups. The trick is to replace a group $G$ possibly with
torsion by its preimage $\widehat{G}$ in the universal central extension of
$\homeo^+(S^1)$, which is a torsion--free group. 
Then $\widehat{G}$ is a central extension of $G$ by some element $T$. We can then
add an auxiliary generator $\beta$ and a relation $\beta T \beta^{-1} = T^2$. This
produces a new torsion--free group $\widehat{G}'$ containing $\widehat{G}$ as
a subgroup. Now,
under any homomorphism from $\widehat{G}'$ to $\homeo^+(S^1)$, the image of
$T$ has rotation
number zero, and therefore has a fixed point. Since $T$ is central in
$\widehat{G}$, the image of $\widehat{G}$ must preserve the fixed point set of $T$.
By blowing this set up, and the complementary regions to the fixed point set of $T$
down, one sees that the action of $\widehat{G}$ on the fixed point set of $T$ is
semi--conjugate to an action which factors through $G$, and therefore has the desired
dynamics.
\end{rmk}

\subsection{Transcendental rotation numbers and arithmetic lattices}

In this section we will introduce another construction which lets us
force certain trancendental rotation numbers. Explicitly, these rotation numbers
will be of the form $\cos^{-1}(r)/2\pi$ where $r$ is a real number of
absolute value $<1$ contained in certain algebraic number fields.

We will give the definition of an {\em arithmetic Fuchsian group}. This
definition contains several terms which may be unfamiliar to the reader,
and therefore we follow the definition with a brief exposition of the
terms involved. Our discussion follows Borel \cite{Borel}
and Vigneras \cite{Vig}. Also, \cite{ReidMac} is an excellent general reference.

\begin{defn}
Let $F$ be a totally real number field and $A$ a quaternion algebra over $F$ which
is ramified at all but one archimedean place. Then a Fuchsian group $\Gamma$ is
{\em arithmetic} if for the embedding $\sigma: A\otimes_F \R \to M_2(\R)$
induced by the unramified archimedean place $F \to \R$,
there is an order $\O$ in $A$ such that $\Gamma$ is commensurable with
$P\sigma(\O^1)$ where $\O^1$ is the group of elements in $\O$ of norm $1$.
\end{defn}

Here a number field is {\em totally real} if all of its Galois embeddings
are subfields of $\R$. A {\em quaternion algebra} $A$ over $F$ is a $4$--dimensional
algebra, with center equal to $F$, and trivial radical. Every such algebra can
be denoted
$$A = \Bigl( \frac{a,b} F \Bigr)$$
where $a,b \in F$, and $A$ is generated (additively) over $F$ by $1,i,j,k$ where
$$i^2 = a, \; j^2 = b, \; k = ij = -ji$$ 
If $x = x_0 + x_1i + x_2j + x_3k$ is an arbitrary element of $A$, then the
{\em trace} of $x$ is $2x_0$, and the {\em norm} is $x_0^2 -x_1^2a - x_2^2b + x_3^2ab$

The {\em archimedean places} of $F$ are the completions $F \to \R$ coming
from the different Galois embeddings of $F$ in $\R$. Such a completion induces an
inclusion of $A$ into a quaternion algebra over $\R$. The only two such algebras,
up to isomorphism, are the ring of $2 \times 2$ real matrices $M_2(\R)$, and the
ring of Hamilton's quaternions $\H$. We say an archimedean place is {\em ramified}
if $A \otimes_F \R \cong \H$ (where $F$ is identified with its image as a subfield of $\R$
under the relevant Galois embedding) and unramified otherwise.

Let $\O_F$ denote the ring of integers in $F$.
An {\em order} $\O$ in a quaternion algebra is a subring of $A$ containing $1$ which
is a finitely generated $\O_F$ module generating the algebra $A$ over $F$. The
elements $\O^1$ of $\O$ of norm $1$ are a group under multiplication.

Parsing the statements above, one sees that for $A,F$ satisfying the conditions of
the definition, the various Galois embeddings $\rho_0,\dots,\rho_n$ of $F$ in
$\R$ define embeddings (which we also denote by $\rho_i$) which, after re--ordering,
we can assume are of the form
$$\rho:A \to \rho_0(A) \times \rho_1(A) \times \dots \times \rho_n(A) \subset
M_2(\R) \times \H \times \dots \times \H$$
and the image of $\O^1$ lies in the product
$$\rho(\O^1) \subset \SL(2,\R) \times \SU(2) \times \dots \times \SU(2)$$
(here we have identified $\SU(2)$ with the group of Hamilton's unit quaternions).
It is not too hard to show that the image of $\O^1$ under $\rho$ is always a
{\em lattice}; since all but one factor of the image is compact, it implies that
the image in $\SL(2,\R)$ is a lattice, and therefore after projection, the image in
$\PSL(2,\R)$ is also a lattice.

\begin{exa}
Let $\Delta'$ be fundamental group of the
hyperbolic triangle orbifold with cone angles $\pi/2,\pi/4,\pi/8$.
Then $\Delta'$ is derived from a quaternion algebra over $\Q(\sqrt{2})$.
The group $\Delta'$ can be conjugated in $\PSL(2,\R)$ into
$\PSL(2,\Q(2^{1/4}))$, and the traces of elements lie in $\Q(\sqrt{2})$.
\end{exa}

\begin{thm}\label{arithmetic_rotation_numbers}
Let $F$ be a totally real number field, and $A$ a quaternion algebra over
$F$ which is ramified at all but one archimedean place $\rho:F \to \R$.
Let $q \in A$ satisfy $\text{norm}(q)=1$ and
$|\text{trace}(\rho(q))| < 2$. Let $$\theta =
\frac{\cos^{-1}(\text{trace}(\rho(q))/2)}{\pi}$$
Then there is a finitely presented group $G_\theta$ and an element
$\alpha \in G_\theta$ such that the set of rotation numbers of $\alpha$
under homomorphisms from $G_\theta$ to $\homeo^+(S^1)$ is exactly the set
$\lbrace 0, \theta, -\theta \rbrace$.
\end{thm}
\begin{proof}
Choose an order $\O$ with group of units $\O^1$. Then $\Gamma = P\rho(\O^1)$
is a lattice of cofinite volume in $\PSL(2,\R)$. Suppose first that
$\Gamma$ is cocompact. Then by the constructions in \S 4.3 we can embed
$\Gamma$ in a larger group containing a copy of $\Delta$, and conjugate the
action of some finite index subgroup of $\Gamma$ to some finite index subgroup of
$\Delta$. This forces $\Gamma$ to act on $S^1$ in a manner semiconjugate to
its tautological representation in $\PSL(2,\R) < \homeo^+(S^1)$, by
theorem~\ref{maximal_semiconjugate}, or else to act trivially.

The element $q \in A$ conjugates the order $\O$ to another order $\O'$ in $A$.
Now, any two orders in a quaternion algebra are commensurable. This
is straightforward to prove; for details consult \cite{Vig}.

But this implies that the projection $P\rho(q) \in \PSL(2,\R)$ is
an element of the commensurator of $\Gamma = P\rho(\O^1)$, and there are
(explicit) finitely presented subgroups $\Gamma_q^1,\Gamma_q^2$ of $\Gamma$
such that $P\rho(q) \Gamma_q^1 P\rho(q)^{-1} = \Gamma_q^2$.
So add a generator $\gamma_q$ and relations
$\gamma_q \Gamma_q^1 \gamma_q^{-1} = \Gamma_q^2$.

Note that the dynamics of $\gamma_q$ on $S^1$ are determined up to semi--conjugacy
by its action on the set of fixed points of elements of
$\Gamma_q^1$. Since such fixed points
must be in the correct circular order for any nontrivial representation of our group
to $\homeo^+(S^1)$, the dynamics of $\gamma_q$ are likewise fixed, up to
semi--conjugacy, and we just need to calculate the rotation number of $\gamma_q$ for
the tautological representation in $\PSL(2,\R) < \homeo^+(S^1)$.

There is a subtle issue, which is that if $\sigma|_\Gamma: \Gamma \to \homeo^+(S^1)$
is the trivial representation, then the dynamics of $\sigma(\gamma_q)$ are {\em not}
constrained by the fact that $\gamma_q \Gamma_q^1 \gamma_q^{-1} = \Gamma_q^2$.
So we must add an extra relation: let $q' \in A$ be such that 
$q^{-1} q' \in \O^1$ and $|\text{trace}(P\rho(q'))|>2$. Then we
add a generator $\gamma_{q'}$, and relations of the form
$\gamma_{q'} \Gamma_{q'}^1 \gamma_{q'}^{-1} = \Gamma_{q'}^2$ and also
a relation $\gamma_q \gamma = \gamma_{q'}$ for some appropriate
$\gamma \in \Gamma$. Add too a new auxiliary generator $\delta$ and the
relation $\delta \gamma_{q'} \delta^{-1} = \gamma$, which will certainly
be satisfied for some $\delta$ for any nontrivial representation of $\Gamma$.
This does the following:
if the representation $\sigma$ is trivial on $\Gamma$, then
$\sigma(\gamma_q) = \sigma(\gamma_{q'}) = \sigma(\gamma) = \id$, as it should
be.

Otherwise, the representation $\sigma|_\Gamma$ is semi--conjugate to the
tautological representation, and the
rotation number of $\sigma(\gamma_q)$ is equal to the rotation number
of $P\rho(q) \in \PSL(2,\R)$, which is 
$\theta = \cos^{-1}(\text{trace}(P\rho(q))/2)/\pi$. It follows that we have
proved the theorem when $\Gamma$ is cocompact. Note that we must divide by
$\pi$ in the denominator rather than $2\pi$, because we are calculating rotation
numbers in $\PSL(2,\R)$ rather than in $\SL(2,\R)$.

Now, if $\Gamma$ is {\em not} cocompact, the orbifold $\H^2/\Gamma$ does not have
a fundamental cycle in $H^2(\Gamma;\R)$ (even virtually), so we cannot directly apply 
the constructions of the previous section. But there is a simple doubling trick:
If $O$ is the orbifold $\H^2/\Gamma$, then there is a homeomorphic hyperbolic
surface $O_c$ where each boundary cusp has been replaced by a geodesic loop.
Then double $O_c$ to $D(O_c)$, and use the constructions from $\S 4.3$ to
embed $\pi_1(D(O_c))$ in a larger group in such a way that $\pi_1(D(O_c))$ is
forced to act semi--conjugately to a Fuchsian action. Then every action
induces an action of the subgroup 
$\Gamma = \pi_1(O_c) <\pi_1(D(O_c))$ on $S^1$ which is semi--conjugate
to the desired action. As before, any $\rho(q)$ conjugating one finite index
subgroup of $\Gamma$ to another must have the correct dynamics on the set of fixed
points of $\Gamma$. It follows that for such a $q$, $\gamma_{q}$ must be
semi--conjugate to its image under the tautological representation.
In particular, the rotation number is forced, and the theorem is proved in this
case too.
\end{proof}

\begin{cor}
Rotation numbers $\theta$ can be forced by finitely presented groups $G_\theta$
for countably many algebraically
independent transcendental numbers.
\end{cor}
\begin{proof}
If $\theta$ as in theorem~\ref{arithmetic_rotation_numbers} is not rational, it
is transcendental, by the celebrated theorem of Gel'fond and Schneider 
(\cite{Gelfond},\cite{Schneider}), which states that for any $\alpha \in \C$
algebraic, either $i\log(\alpha)/\pi$ is in $\Q$, or it is transcendental.
For $\alpha = \cos(2\pi\theta) + i\sin(2\pi\theta)$, $\alpha$ is algebraic iff
$\alpha + \alpha^{-1} = 2\cos(2\pi\theta)$ is algebraic, and then 
$i\log(\alpha)/\pi = -2\theta$ is
either rational or transcendental.

The other statement
follows from the existence of infinitely many distinct totally real number fields $F$
over which one can define arithmetic Fuchsian groups. See \cite{Vig}.
\end{proof}

\begin{rmk}
For $F$ a nontrivial totally real
extension of $\Q$, any arithmetic lattice in $\PSL(2,\R)$ defined over $F$ is
cocompact. So the only non--cocompact lattices to consider in the proof of 
theorem~\ref{arithmetic_rotation_numbers} are those commensurable with $\PSL(2,\Z)$.
\end{rmk}

\begin{rmk}
An alternative construction of rigid subgroups of $\homeo^+(S^1)$ with
many symmetries comes from \cite{Ghys_lattice}. Ghys shows that if $\Gamma$ is
an irreducible lattice in a semi--simple Lie group of real rank at least $2$
of the form
$\SL(2,\R)^k \times G$ where $G$ has no $\SL(2,\R)$ factor,
then any nonfinite action of $\Gamma$ on $S^1$ is semi--conjugate to an action
obtained by projecting onto one of the $\SL(2,\R)$ factors.

Our construction shows that rank $1$ lattices can also be ``rigidified'' by
embedding them in slightly larger groups, in a straightforward way.
\end{rmk}

\subsection{Arithmetic and $l$--arithmetic of rotation numbers}

In this subsection we introduce some tools for producing new rotation
numbers from the rotation numbers we produced in \S 4.3 and \S 4.4.

\begin{lem}\label{add_rotation_numbers}
Let $\theta_1,\theta_2$ be rotation numbers forced by $G_{\theta_1}, G_{\theta_2}$
as in theorem~\ref{rational_rotation_numbers} and
theorem~\ref{arithmetic_rotation_numbers}. Then $\theta_1 + \theta_2$ can be
forced.
\end{lem}
\begin{proof}
If $\alpha_i \in G_{\theta_i}$ are the elements whose rotation numbers are
forced to be equal to one of $\lbrace 0, \pm \theta_i \rbrace$ respectively,
then we can form the group
$$G = \langle G_{\theta_1}, G_{\theta_2} \; | \; [\alpha_1,\alpha_2] = \id \rangle$$
If two elements $\mu,\nu \in \homeo^+(S^1)$ commute, then it is easy to see
that the rotation number of the product satisfies 
$$\rot(\mu \nu) = \rot(\mu) + \rot(\nu)$$
$G$ is not quite the group we want, since the element $\alpha_1\alpha_2$ can
have rotation numbers 
$\lbrace 0,\pm\theta_1,\pm\theta_2,\pm \theta_1 \pm \theta_2 \rbrace$ under
various representations.

However, the groups $G_{\theta_1}$ and $G_{\theta_2}$ contain canonical copies of
the subgroup $\Delta$. Call these copies $\Delta_1,\Delta_2$. If we add a
generator $\gamma$ and a relation $\gamma \Delta_1 \gamma^{-1} = \Delta_2$
then the rotation number of $\alpha_1$ vanishes iff the rotation number of
$\alpha_2$ vanishes, and their signs are similarly correlated.
It follows that 
$$G_{\theta_1 + \theta_2} = \langle G,\gamma \; | \; \gamma \Delta_1 \gamma^{-1}
= \Delta_2 \rangle$$
forces the element $\alpha_1 \alpha_2$ to have the relevant dynamics.
\end{proof}

Similarly, rotation numbers can be subtracted, by forcing conjugacy between
oppositely oriented copies of an amphichiral index $2$ subgroup of $\Delta$.
It follows that the set $X$ of rotation numbers that we can force by our
techniques is a $\Z$--module.

However, we can also introduce a different kind of ``addition'' of rotation
numbers, dependent on a parameter, which can be thought of as a deformation of
usual addition.

\begin{defn}
Let $\theta_1,\theta_2 \in \R/\Z$ and $l \in \R$. Define
$$\theta_1 +_l \theta_2 = \frac {\cos^{-1} (\cos(\pi\theta_1)\cos(\pi\theta_2) 
- \cosh(l)\sin(\pi\theta_1)\sin(\pi\theta_2))} {\pi}$$
\end{defn}

The formula is a bit clumsy, since we have normalized rotation numbers to
live in $\R/\Z$ instead of the more natural $\R/2\pi\Z$.
Note that $+_0 = +$, and that $+_l = +_{-l}$
so this formula really should be thought of as a symmetric deformation
of usual addition. Notice too that for a given $l$, $+_l$ is only
defined for sufficiently small $\theta_1,\theta_2$, but that for any $l\ne 0,\theta_1$,
it is defined for $\theta_2$ in the complement of some neighborhood of $-\theta_1$.

The geometric meaning of $+_l$ is simple: if $\alpha_1,\alpha_2$ are
rotations of $\H^2$ through angles $2\pi\theta_1,2\pi\theta_2$
about points $p_1,p_2$ which are a distance $l$ apart,
then their product is a rotation through angle $2\pi(\theta_1 +_l \theta_2)$
for suitable values of $\theta_1,\theta_2,l$.

By abuse of notation, we call a set $X \subset \R/\Z$ which is
closed with respect to the (partially
defined) operation $+_l$ for some $l$ a {\em $+_l$--module}.
If $L$ is a subset of $\R$, we refer to a set closed with respect to
$+_l$ (where defined) for every $l \in L$ as an {\em $+_L$--module}.

\begin{thm}\label{deformed_addition}
Let $Y$ be the set of rotation numbers forced by 
theorem~\ref{rational_rotation_numbers} and theorem~\ref{arithmetic_rotation_numbers}.
Then there is a dense set $L$ containing $0$ of real numbers $l \in \R$ which are of the form
$\log(r)$ for $r$ algebraic, so that every element of the $+_L$--module $X$
generated by $Y$ can be forced.
\end{thm}
\begin{proof}
Let $A,F,\O^1,\Gamma$ be as in theorem~\ref{arithmetic_rotation_numbers}.
Then the commensurator of $\Gamma$ is dense in $\PSL(2,\R)$. Let
$C(\Gamma) \subset \PSL(2,\R)$ denote the commensurator of $\Gamma$.
Then there exist pairs of elliptic elements $q_1,q_2$ with
irrational rotation number which are conjugate
in $C(\Gamma)$ whose centers are distance $r$ away, for a dense set of $l \in \R$.
Let $q \in C(\Gamma)$ conjugate $q_1$ to $q_2$.
The elements $q_1,q_2,q$ can be characterized up to semi--conjugacy by
the fact that they conjugate specific
finite index subgroups of $\Gamma$ to each other. So we can produce a finitely
presented group $\Gamma(q,q_1,q_2)$
which contains a copy of $\Gamma$, and elements $q_1,q_2,q$
which must act on $S^1$, up to semi--conjugacy, in the desired manner.
Note that the only elements in $\homeo^+(S^1)$ which commute with an
{\em irrational} rigid rotation are other rigid rotations. This follows
from the well--known fact that Lebesgue measure on $S^1$ is 
the unique invariant measure for an irrational rigid rotation. In particular,
any element of $\homeo^+(S^1)$ which commutes with some $q_i$ must be semi--conjugate
to a rigid rotation with respect to the hyperbolic visual measure on $S^1$ centered
at the center of $q_i$.

If $G_{\theta_1},G_{\theta_2}$ have been produced, containing canonical
subgroups isomorphic to $\Delta$, and elements $\alpha_1,\alpha_2$ whose
rotation numbers are forced to be equal to $\theta_1,\theta_2$ respectively
(after fixing a nontrivial homomorphism $\Delta \to \homeo^+(S^1)$),
then we can introduce relations saying that $\alpha_i$ and $q_i$ commute,
for $i=1,2$.

Define a group $G$, generated by $\Gamma(q,q_1,q_2),G_{\theta_1},G_{\theta_2}$
and some auxiliary generators which conjugate the canonical copies of
$\Delta$ in each group to each other. Add relations $[\alpha_i,q_i] = \id$ for
$i=1,2$. Then every action of $G$ on $S^1$ which restricts to a nontrivial
action of $\Delta$ is semi--conjugate to an action
in which $\alpha_1,\alpha_2$ are rotations centered at the centers of
$q_1,q_2$ under the tautological representation to $\PSL(2,\R)$. In particular,
the composition $\alpha_1\alpha_2$ has rotation number $\theta_1 +_l \theta_2$,
as required.
\end{proof}

The operation $+_l$ might seem exotic, and the fact that it is only partially
defined might seem like a detriment. But we will see
that this apparent deficit is actually an asset: the operation $+_l$
lets us impose geometric {\em inequalities}
on rotation numbers. This is important, because it lets us {\em divide}
rotation numbers, and gives the
set of forceable rotation numbers the structure of a $\Q$--module.

This requires some explanation: if $\theta \in \R/\Z$, then for any integer
$p$, $\theta/p$ is only well--defined up to multiples of $1/p$. So to say
that the set $X$ of forceable rotation numbers has the structure of a $\Q$--module
means that for any $\theta \in X$, and any $p \in \N$, that all $p$ possible
values of $\theta/p$ are in $X$.

Obviously, if we have $\alpha \in G$ such that the set of rotation numbers
$\rho(\alpha)$ as $\rho$ varies in $\hom(G;\homeo^+(S^1))$ is equal to the set
$\lbrace 0,\pm \theta \rbrace$, then if we define $G^{1/p}$ to be the group
$$G^{1/p} = \langle G,\beta \; | \; \beta^p = \alpha \rangle$$
then the set of rotation numbers $\rho(\beta)$ as $\rho$ varies in
$\hom(G^{1/p};\homeo^+(S^1))$ is equal to the set
$$\lbrace 0, \pm \theta/p, \pm (\theta+1)/p, \dots \pm (\theta + p - 1)/p \rbrace$$
If we can further constrain the rotation number of $\beta$ to satisfy certain
appropriately chosen {\em inequalities}, we can force the rotation number $q\theta/p$ that
we want.

For a pair of distinct numbers $t_1,t_2 \in \R/\Z$,
we let $[t_1,t_2]$ denote the oriented interval from $t_1$ to $t_2$.

\begin{lem}\label{force_inequality}
There is a dense subset $T \subset \R/\Z\times \R/\Z$ such that for any $t \in T$
with co--ordinates $t_1,t_2$, there is a pair $G_{[t_1,t_2]},\alpha$ 
where $G_{[t_1,t_2]}$ is finitely generated and
$\alpha \in G_{[t_1,t_2]}$, such that the set of rotation numbers $\rho(\alpha)$
as $\rho$ varies over elements of
$\hom(G_{[t_1,t_2]},\homeo^+(S^1))$ is exactly the set
$\lbrace 0\rbrace \cup \pm[t_1,t_2]$.
\end{lem}
\begin{proof}
For a given $r \in \R$ and $q \in \R/\Z$, the set of $q' \in \R/\Z$ such
that $q +_r q'$ is well--defined and non--zero
is a connected open interval $I_{r,q}$ whose complement $\R/\Z - I_{r,q}$ 
contains $-q$. 
Moreover, the size and position of this interval varies continuously
as a function of $r,q$ and can be made as large or small as desired, and
placed arbitrarily in $\R/\Z$, by
suitable choice of $r,q$.

As before, let $A,F,\O^1,\Gamma$ be as in theorem~\ref{arithmetic_rotation_numbers}.
Then the commensurator $C(\Gamma)$
of $\Gamma$ is dense in $\PSL(2,\R)$, and we can choose
elements $\alpha,\alpha',\beta \in C(\Gamma)$ such that $\alpha,\alpha'$ are
rotations through angle $q+\epsilon$
which are conjugated by $\beta$, which has translation length $r+\epsilon$
for $\epsilon$ as small as desired. Then if we add a generator $\gamma$ and
a relation $[\gamma,\alpha'] = \id$, and another relation which forces
$\rot(\alpha\gamma)=0$, then the rotation number of $\gamma$
must be in $\R/\Z - I_{r,q}$, and
every number in this interval is realized by some $\gamma$.
Here the relation $\rot(\alpha\gamma)=0$ can be forced by adding another generator
$\mu$ and a relation $\mu \alpha \gamma \mu^{-1} = \beta$, as in the proof of
theorem~\ref{arithmetic_rotation_numbers}.

It follows that we have constructed a finitely generated group $G$ and an
element $\gamma \in G$ such that the set of rotation numbers of $\gamma$ is
exactly $\lbrace 0 \rbrace \cup \pm(\R/\Z - I_{r,q})$, and if $0 \in I_{r,q}$,
the rotation number $0$ is attained iff the representation restricts to
the trivial representation on the canonical $\Delta$ subgroup of $G$. This proves the lemma.
\end{proof}

This lemma immediately implies that the set of forceable
numbers we can construct by a combination of our methods is a $\Q$--module
(in the sense above). Moreover, this construction can be generalized in
the following way.

We can think of the operation of {\em dividing by $2$} as solving an
algebraic equation with respect to the operation $+$. That is,
$x = y/2$ can be re--written as $x+x=y$, which can be solved for
$x$, and the various (discrete) solutions discriminated amongst by
applications of lemma~\ref{force_inequality}. We can also force
solutions to more general equations with coefficients in the
set of rotation numbers we have already forced, and with operations
$+_l$ with $l$ as above. As long as the set of solutions is discrete,
we can discriminate amongst them by lemma~\ref{force_inequality}.

In this way, we define the {\em algebraic closure} of a set of
rotation numbers with respect to the set of operations $+_l$.

\begin{defn}
If $Y \subset S^1$ is a set of rotation numbers and $L \subset \R$, 
the {\em algebraic closure of $Y$} with respect to the operations $+_l$
with $l \in L$, is defined inductively as
the smallest set which includes $Y$ itself, and also includes every
solution to every finite equation or system of equations
(with discrete solution set) in a variable
$x$, and possibly auxiliary variables $y_1,\dots,y_n$,
and operations $+_l$, and coefficients in the algebraic closure of
$Y$.
\end{defn}

The reason to allow auxiliary variables $y_i$ is that there does not
seem to be a simple way to eliminate variables with respect to the operations $+_l$.

\begin{exa}
The following are examples of equations and systems of equations
to be solved in $x$:
\begin{enumerate}
\item{$x+_l x = \theta$}
\item{$x+_{l_1} x = \theta +_{l_2} x$}
\item{$(y +_{l_1} ((x +_{l_2} x) +_{l_3} \theta_1) = \theta_2 +_{l_4} y, \; 
x +_{l_5} y = y +_{l_6} \theta_3$}
\end{enumerate}
\end{exa}

\begin{transcendental_module_cor}
Let $X$ be the set of forceable rotation numbers. Then $X$ contains
countably infinitely many algebraically independent transcendental
numbers, as well as all the rational numbers. Moreover,
there is a dense set $L$ of real numbers $l \in \R$ containing $0$,
which are of the form
$\log(r)$ for $r$ algebraic, so that $X$ is algebraically closed
with respect to the operations $+_l$.
\end{transcendental_module_cor}
\begin{proof}
For the set of rotation numbers that we have forced so far, this is
an easy consequence of lemma~\ref{force_inequality}. If $(G_i,\alpha_i)$ are some
collection of pairs which force rotation numbers $\theta_i$ respectively,
we need to control the fact that the rotation numbers $\theta_i$ are only
really forced {\em up to sign}, and therefore we must distinguish between
numbers $\theta_i + \theta_j, \theta_i - \theta_j$ and so on.
But this can also easily be accomplished by using lemma~\ref{force_inequality}
to correlate the signs of the rotation numbers of the $\alpha_i$. Having made this
observation, the corollary follows easily from the constructions above.
\end{proof}

The following theorem can be thought of as showing that 
{\em finite unions of forceable sets of rotation numbers are forceable}:

\vfill
\pagebreak

\begin{outer_approximation_thm}[Outer approximation theorem]\label{outer_approximation}
Let  $K$ be any closed subset of $\R/\Z$. Then there are
a sequence of pairs $G_{K_i},\alpha_i$ where $G_{K_i}$ is
a finitely presented group and $\alpha_i \in G_{K_i}$, and
closed subsets $K_i$ of $\R/\Z$ such that
\begin{enumerate}
\item{Each $K_{i+1} \subset K_i$.}
\item{The intersection $\cap_i K_i = K$.}
\item{The set of rotation numbers $\rho(\alpha_i)$ as $\rho$ ranges over
$\hom(G_{K_i},\homeo^+(S^1))$ is exactly equal to 
$\lbrace 0 \rbrace \cup K_i \cup - K_i$.}
\item{There is a canonical element
$\nu_i \in G_{K_i}$ of order $3$, so that if $\rot(\rho(\nu_i))=0$, $\rho(\alpha_i)$
is trivial, and if $\rot(\rho(\nu_i))=1/3$, the set
of compatible $\rot(\rho(\alpha_i))$ is exactly equal to $K_i$.}
\end{enumerate}
\end{outer_approximation_thm}
\begin{proof}
By lemma~\ref{force_inequality} we can produce $G_I$ for a dense set
of closed intervals
$I \subset \R/\Z$. If $I,J$ are two such intervals, we can produce
$G_{I \cup J}$ as follows. Without loss of generality, we can assume $I \cup J$
is not the entire circle $\R/\Z$, otherwise we our done.

Let $\alpha_I,\alpha_J$ be the distinguished elements
of $G_I,G_J$ whose rotation number has been forced. Let $\Delta_I,\Delta_J$ be
the distinguished subgroups isomorphic to $\Delta$. Let $\mu_I,\mu_J$ be
the distinguished elements in $\Delta_I,\Delta_J$ respectively of order $7$.
Let $\Gamma$ be some fixed arithmetic group, with its own copy $\Delta_\Gamma$
of $\Delta$ and distinguished element $\mu_\Gamma$ of order $7$.

As a first approximation, let $G_{I \cup J}$ be generated by $G_I,G_J$ and $\Gamma$. 
Add relations that $[\mu_I,\mu_J]=\id$ and $\mu_I\mu_J = \mu_\Gamma$. Then add
another generator $\delta$ and a relation
$[\delta \alpha_I \delta^{-1},\alpha_J] = \id$, and
let $\alpha = \delta \alpha_I \delta^{-1} \alpha_J$ be the
product.

Now, if $\rho:G_{I \cup J} \to \homeo^+(S^1)$ is any representation, it
must be trivial on at least one of the copies $\Delta_I,\Delta_J,\Delta_\Gamma$,
since the sum of the rotation numbers satisfies
$$\rot(\rho(\mu_I)) + \rot(\rho(\mu_J)) = \rot(\rho(\mu_\Gamma))$$
We want to force $\rho$ to be {\em nontrivial} on $\Delta_\Gamma$, or else
{\em necessarily trivial} on all three.

Add as generators elements $\iota, \kappa \in C(\Gamma)$ whose
images under the tautological representation of $C(\Gamma)$ are
elliptic and hyperbolic elements of $\PSL(2,\R)$ respectively, where the
rotation number of $\iota$ is $r$, and the center of the
rotation $\iota$ is translated hyperbolic distance $l$ by $\kappa$.
Add another auxiliary generator $\phi$ and a relation
$$[\phi\alpha\phi^{-1},\kappa \iota \kappa^{-1}] = \id$$
Now force the composition $\iota \phi\alpha\phi^{-1}$ to have rotation number $0$ by
adding an auxiliary generator which conjugates $\iota \phi \alpha \phi^{-1}$ to
some hyperbolic element of $\Gamma$.

If $\rho$ is trivial on $\Delta_\Gamma$ then $\phi\alpha\phi^{-1}$, and therefore
$\alpha$, is trivial. Otherwise, $\rho$ is nontrivial on
$\Delta_\Gamma$, and therefore is nontrivial on at most one of
$\Delta_I,\Delta_J$. It follows that, after fixing a nontrivial homomorphism from
$\Gamma$ to $\homeo^+(S^1)$, the set of compatible rotation numbers of
$\alpha$ in $\R/\Z$ is exactly equal to $(I \cup J) \cap (\R/\Z - I_{r,q})$.
By choosing $r,q$ appropriately, we can make $I_{r,q}$ disjoint from
$I \cup J$.

By an obvious induction, we can construct $G_K$ for 
$K$ any finite union of disjoint closed
intervals $I$ where $G_I$ can be forced. This proves the theorem.
\end{proof}

\subsection{The representation topology}

We now have the tools to completely describe the representation topology,
defined in \S 4.2.

\begin{lem}\label{generate_closed_sets}
Let $K$ be a closed subset of $\R/\Z$ which is invariant under $\theta \to -\theta$
and contains $0$. Then there is a countable group $G_K$ and an element
$\alpha \in G_K$ such that the set of rotation numbers of $\alpha$ induced
by actions of $G_K$ on $S^1$ is exactly equal to $K$. Moreover, if $K$ is
algorithmically constructible, so is $G_K$.
\end{lem}
\begin{proof}
We can just take $G_K$ to be the amalgamation of the groups $G_{K_i}$
constructed in theorem~\ref{outer_approximation} over the cyclic subgroups
generated by each $\alpha_i$, which are identified by the isomorphisms taking
$\alpha_i$ to $\alpha_j$.
\end{proof}

\begin{lem}\label{closed_sets_are_closed}
Let $G$ be a countable group, and $\alpha \in G$ some element. Then
the set of rotation numbers of $\alpha$ induced by actions of $G$ on $S^1$
is closed (in the usual sense) as a subset of $\R/\Z$, is invariant under
$\theta \to -\theta$, and contains $0$.
\end{lem}
\begin{proof}
The only assertion requiring proof is that the set of rotation numbers
achieved by $\alpha$ is {\em closed}. But this will follow easily from the
{\em compactness} of the space of semi--conjugacy classes of group
actions, as characterized by Ghys' theorem. Explicitly, let $\rho_i$
be a sequence of representations $\rho_i:G \to \homeo^+(S^1)$ such
that the rotation numbers $\rot(\rho_i(\alpha)) = r_i$ where $r_i \to r \in \R/\Z$,
with convergence in the usual sense. Then the homomorphisms $\rho_i$ determine
bounded cocycles $c_i: G^3/G \to \Z$ taking the values $0,1$.
Since the product $\lbrace 0,1 \rbrace^G$ is compact, we can find a
subsequence which is eventually constant on any finite subset of $G^3/G$,
which therefore converges to some bounded cocycle $c$. By Ghys' 
theorem~\ref{semiconjugacy_invariant}, $c$ determines a representation
$\rho:G \to \homeo^+(S^1)$ up to semiconjugacy. The rotation number
$\rot(\rho(\alpha))$ is completely determined by the restriction of $c$
to the subgroup $\langle \alpha \rangle$ generated by $\alpha$. By construction,
it is clear that it is equal to $r$.
\end{proof}

It follows that we can completely characterize the topology on $S^1$
generated by forceable sets of rotation numbers:

\begin{representation_topology_cor}
The set of subsets of $S^1$ of the form $\rot(X(G,\alpha))$ 
where $G$ varies over all countable groups, and $\alpha \in G$ is
arbitrary, are precisely the nonempty closed subsets of a topology, called the
{\em representation topology}.

The nonempty closed subsets in the representation topology on $S^1$ are
exactly unions $\lbrace 0 \rbrace \cup K$ where $K$ is closed (in the
usual sense) and invariant under $x \to -x$.
\end{representation_topology_cor}

\subsection{Smoothness issues}

We point out that the question of which rotation numbers can be forced
depends significantly on how much smoothness is assumed. The constructions
throughout \S 4 cannot generally be made $C^3$.
This is because we frequently conjugate the action
of various Fuchsian groups (which are not conjugate in $\PSL(2,\R)$)
to each other.

However, we have the following
theorem of Ghys from \cite{Ghys_rigid}:

\begin{thm}[Ghys]
A $C^3$ action of a negatively curved closed surface group $\Gamma$ on $\homeo^+(S^1)$
of maximal Euler class is $C^3$ conjugate to a Fuchsian action. Moreover,
two discrete cocompact subgroups of $\PSL(2,\R)$ are conjugate by a 
$C^1$ diffeomorphism of $S^1$ only when they are conjugate in $\PSL(2,\R)$
\end{thm}

\begin{rmk}
In general, the Gromov boundary of a $-\frac 1 4 < k \le -1$--pinched
negatively curved manifold has a natural $C^1$ structure, but not $C^2$.
See \cite{Sull} for more details.
\end{rmk}

In fact, where we are occasionally forced to construct actions which are
necessarily semi--conjugate (but not conjugate) to Fuchsian actions, one can
show the actions cannot in general even be made $C^2$. However, they {\em can}
generally by made $C^1$, by arguments of Denjoy and Pixton.

It would be very interesting to see how the representation topology on $S^1$
stiffens and rigidifies as the degree of smoothness increases. A recent paper of
Amie Wilkinson and Lizzie Burslem \cite{Wilk_Burs} produces striking examples
of rigidity of solvable group actions at every degree $r$ of smoothness.

\end{document}